\documentclass{article}

\usepackage{arxiv}

\usepackage[utf8]{inputenc} 
\usepackage[T1]{fontenc}    
\usepackage{hyperref}       
\usepackage{url}            
\usepackage{amsfonts}       
\usepackage{lipsum}	
\usepackage{amsmath}
\usepackage{amsthm}

\newtheorem{thm}{Theorem}[section]
 
 \newtheorem{rem}{Remark}[section]
 \newtheorem{lem}{Lemma}[section]

 \numberwithin{equation}{section}

\title{	Existence and regularity  of positive solutions for Schrödinger-Maxwell system with singularity}

\author{
Abdelaaziz Sbai, Youssef El hadfi and Mounim El Ouardy \\
Laboratory LIPIM\\
National School of Applied Sciences Khouribga\\
Sultan Moulay Slimane University, Morocco\\
 \texttt{sbaiabdlaaziz@gmail.com}\\
  \texttt{yelhadfi@gmail.com}\\
   \texttt{ monim.elrdy@gmail.com}\\
}

\begin{document}
\maketitle
\begin{abstract}
	In this paper we are going to prove existence  for positive
	solutions of  the following Schrödinger–Maxwell system of singular elliptic equations:
	\begin{equation}
	\left\{\begin{array}{l}
	u \in W_{0}^{1,2}(\Omega):-\operatorname{div}\left(a(x) \nabla u\right)+\psi|u|^{r-2} u=\frac{f(x)}{u^{\theta}}, \\
	\psi \in W_{0}^{1,2}(\Omega):-\operatorname{div}(M(x) \nabla \psi)=|u|^{r}
	\end{array}\right.
	\end{equation}
	where $\Omega$ is a bounded open set of $\mathbb{R}^{N}, N>2,$ $r>,1,$ $u>0,$ $\psi>0,$  $0 < \theta<1$ and $f$ belongs
	to a suitable Lebesgue space. In particular, we take advantage of the coupling between the two equations of the system
	by demonstrating how the structure of the system gives rise to a regularizing effect on the summability of the solutions.
\end{abstract}
\section{Introduction}
In this paper, we consider the following Schrödinger-Maxwell system with singular term:	
\begin{equation}\label{133}
\left\{\begin{array}{l}
u \in W_{0}^{1,2}(\Omega),u>0 \,:-\operatorname{div}\left(a(x) \nabla u\right)+\psi|u|^{r-2} u=\frac{f(x)}{u^{\theta}}, \\
\psi \in W_{0}^{1,2}(\Omega),\psi>0 \,:-\operatorname{div}(M(x) \nabla \psi)=|u|^{r}
\end{array}\right.
\end{equation}
We will suppose that $\Omega$ is a bounded open set of $\mathbb{R}^{N}, N>2,$ that $r>1$ and that $f$ belongs to $L^{m}(\Omega),$ for some $m > 1$, $0< \theta<1$. Furthermore, the function $a: \Omega \rightarrow \mathbb{R}$ will be a measurable function, such that there exist $0<\alpha \leq \beta$ such that:
\begin{equation}\label{134}
0<\alpha \leq a(x) \leq \beta \,\,\,\,	\mbox{almost everywhere in}\,\,\, \Omega,
\end{equation}
while $M: \Omega \rightarrow \mathbb{R}^{N^{2}}$ will be a measurable matrix, such that:
\begin{equation}\label{135}
M(x) \xi \cdot \xi \geq \alpha|\xi|^{2}, \quad|M(x)| \leq \beta,
\end{equation}
for almost every $x$ in $\Omega,$ and for every $\xi$ in $\mathbb{R}^{N}$. Let us briefly recall the mathematical framework concerning problem \eqref{133}.\\
In the last few decades,
the existence and  regularity of positive solutions to the singular elliptic equation with  singular term $ s^{-\theta},$ $(\theta\in ( 0, 1))$ a have been widely investigated by many researchers, and it seems almost impossible for us to give a complete list of references. We refer the readers to Refs.\cite{b1,b2,s1,s2}.\\
In the case $\theta=0$ many works have appeared concerning the existence and regularity of elliptic systems.
Boccardo in \cite{09} has been studied the existence and regularity results of of elliptic systems  problem;
\begin{equation}\label{i12}
\left\{\begin{array}{ll}
-\Delta u+A \varphi|u|^{r-2} u=f, & u \in W_{0}^{1,2}(\Omega), \\
-\Delta \varphi=|u|^{r}, & \varphi \in W_{0}^{1,2}(\Omega),
\end{array}\right.
\end{equation}
where $\Omega$ is an open bounded subset of $\mathbb{R}^{N}$ with $N>2, A>0$ and $r>1$.
As for the \eqref{i12}, the authors studied the  existence of a weak solution $(u, \varphi)$ in $W_{0}^{1,2}(\Omega) \times W_{0}^{1,2}(\Omega)$ is proved if $f$ belongs to $L^{m}(\Omega)$, with $m \geq \frac{2 N}{N+2}=\left(2^{*}\right)^{\prime}$, where $2^{*}$ is the Sobolev exponent, using once again that $(u, \varphi)$ is a critical point of a suitable functional. The author proves that if $\left(2^{*}\right)^{\prime} \leq m<\frac{2 N r}{N+2+4 r}$, with $r>2^{*}-1$, the second equation of \eqref{i12} admits finite energy solutions even if the datum $|u|^{r}$ does not belong to the dual space $L^{\frac{2 N}{N+2}}(\Omega)$. \\
As for system \eqref{i12},  the solutions $u$ and $ \varphi$ given  by \cite{09} can be seen, if $m \geq \frac{2 N}{N+2}$, as a critical point of saddle type for the indefinite functional:
$$J(u, \varphi)=\frac{1}{2} \int_{\Omega} M(x) \nabla u \nabla v-\frac{A}{2 r} \int_{\Omega} M(x) \nabla \varphi \nabla \varphi+\frac{A}{r} \int_{\Omega} \varphi^{+}|u|^{r}-\int_{\Omega} f u,$$
defined for those functions $u$ and $\varphi$ in $W_{0}^{1,2}(\Omega)$, such that $\varphi^{+}|u|^{r}$ belongs to $L^{1}(\Omega)$ (and $+\infty$ otherwise).
On the other hand, the authors proved in \cite{09} the existence of solutions for the following nonlinear elliptic system
that generalizes \eqref{i12}
\begin{equation}
\left\{\begin{array}{ll}
-\operatorname{div}\left(|\nabla u|^{p-2} \nabla u\right)+A \varphi|u|^{r-2} u=f, & u \in W_{0}^{1, p}(\Omega) \\
-\operatorname{div}\left(|\nabla \varphi|^{p-2} \nabla \varphi\right)=|u|^{r}, & \varphi \in W_{0}^{1, p}(\Omega)
\end{array}\right.
\end{equation}
where $\Omega$ is an open bounded subset of $\mathbb{R}^{N}$ with $N \geq 2,1<p<N, A>0, r>1.$\\
Inspired by the above articles, the main novelty in the present work is to show that the term $\frac{1}{u^{\theta}}$ has a "regularzing effect"
in the sense that the problem \eqref{133} has a distributional solution for all $f\in L^{m}(\Omega)$ with $m>1.$ This term provokes some mathematical difficulties, which make the study of system
\eqref{133} particularly interesting. To our knowledge, the  Schrödinger–Maxwell system with  singular term has not been studied. \\
The plan of the paper is as follows: in the next section we will study a problem
approximating \eqref{133}, proving existence of solutions. Once again by using the approximation scheme we prove estimates that allow us to pass to the limit in the approximate equations and to prove the existence of a weak solution of \eqref{133}.\\
In  the last part of this Section, we prove the existence of a saddle point $(u,\varphi)$ of the following
functional
\begin{equation}\label{1.6}
J(u, \varphi)=\left\{\begin{aligned}
\frac{1}{2} \int_{\Omega} a(x)|\nabla u|^{2}-\frac{1}{2 r} \int_{\Omega} M(x) \nabla \varphi \nabla \varphi & \text { if } \varphi^{+}|u|^{r} \in L^{1}(\Omega), \\
+\frac{1}{r} \int_{\Omega} \varphi^{+}|u|^{r}-\frac{1}{1-\theta}\int_{\Omega} f (u^{+})^{1-\theta}, & \\
+\infty, & \text { otherwise. }
\end{aligned}\right.
\end{equation}
defined on $W_{0}^{1,2}(\Omega)\times W_{0}^{1,2}(\Omega).$
Finally, in the Appendix, we give the proof of an existence result for the first equation of
approximating system and some results allowing to prove the existence of system \eqref{133}.\\
\textbf{Notations.} For a given function $v$ we denote by $v^{+}=\max (v, 0)$ and by $v^{-}=-\min (v, 0) .$  For a fixed $k>0,$ we define the truncation functions $T_{k}: \mathbb{R} \rightarrow \mathbb{R}$ and $G_{k}: \mathbb{R} \rightarrow \mathbb{R}$ as follows
$$
\begin{array}{l}
T_{k}(s):=\max (-k, \min (s, k)) \\
G_{k}(s):=(|s|-k)^{+} \operatorname{sign}(s)
\end{array}
$$
we will also make use of the notation
$$
\int_{\Omega} f(x) d x=\int_{\Omega} f
$$
If no otherwise specified, we will denote by C serval constants whose value may
change from line to line and, sometimes, on the same line. These values will only depend on the data (for instance C can depend on $\Omega$,$\theta$,N,k,...) but the
will never depend on the indexes of the sequences we will often introduce.

\section{ A priori estimates and main results}

In this section we are interesting to prove regularity of u solution of \eqref{133} when the datum f belong to $L^{m}(\Omega),$ with $m> 1$	
\begin{thm}\label{th21}
	Let $0< \theta <1$ and let a and $M$ be such that \eqref{134} and \eqref{135} hold. Let $r>1$ and let $f$ in $L^{m}(\Omega) .$ We have the following:\\
	(i) if $r \geq \frac{2 N}{\theta(N-2)+N+2}$, and if $m \geq\big(\frac{r+1}{1-\theta}\big)^{\prime}$, there exist $u$ and $\psi$ in $W_{0}^{1,2}(\Omega)$, solutions of \eqref{133}; furthermore,\\
	(a) if $m>\frac{N}{2}$, then u belongs to $L^{\infty}(\Omega)$;\\
	(b) if $\big(\frac{r+1}{1-\theta}\big)^{\prime} \leq m<\frac{N}{2}$, then $u$ belongs to $L^{\sigma}(\Omega)$, with $\sigma=$ $\max \left((1+\theta) m^{* *}, \frac{m(2 r+1+\theta)}{m+1} \right)$, where $m^{**}=\frac{N m}{N-2 m}.$\\
	(ii) if $1<r<\frac{2 N}{\theta(N-2)+N+2}$, and if $m \geq \frac{2 N}{\theta(N-2)+N+2}$, there exist $u$ and $\psi$ in $W_{0}^{1,2}(\Omega)$, solutions of \eqref{133} ; \\furthermore:\\
	(a) if $m>\frac{N}{2}$, then $u$ belongs to $L^{\infty}(\Omega)$;\\
	(b) if $\frac{2 N}{\theta(N-2)+N+2} \leq m<\frac{N}{2}$, then $u$ belongs to $L^{(1+\theta)m^{* *}}(\Omega)$, where $m^{* *}=\frac{N m}{N-2 m}$.
\end{thm}
\subsection{The Approximated Problem}
Let $n$ in $\mathbb{N},$ and let $f_{n}=T_{n}(f),$ so that $\left\{f_{n}\right\}$ is a sequence of $L^{\infty}(\Omega)$ functions, which strongly converges to $f$ in $L^{m}(\Omega),$ and satisfies the inequality $\left|f_{n}\right| \leq|f| .$ Thanks to Theorem \ref{3.4} (see the Appendix), for every $n$ in $\mathbb{N},$ there exist weak solutions $u_{n}$ and $\psi_{n}$ in $W_{0}^{1,2}(\Omega) \cap L^{\infty}(\Omega)$, with $\psi_{n} \geq 0,$ of the approximate system:
\begin{equation}\label{eq1}
\left\{\begin{array}{l}
u_{n} \in W_{0}^{1,2}(\Omega):-\operatorname{div}\left(a(x) \nabla u_{n}\right)+\psi_{n}\left|u_{n}\right|^{r-2} u_{n}=\frac{f_{n}}{(\frac{1}{n}+u_{n})^{\theta}},\,\,\,\,\,\, (I) \\
\psi_{n} \in W_{0}^{1,2}(\Omega):-\operatorname{div}\left(M(x) \nabla \psi_{n}\right)=\left|u_{n}\right|^{r},\,\,\,\, (II).
\end{array}\right.
\end{equation}
\subsection{A Priori Estimates}
We are now going to prove some a priori estimates on the sequence of approximated solutions $u_{n}.$
\begin{lem}\label{lc1}
	Let $k>0$ be fixed. The sequence $\left\{T_{k}\left(u_{n}\right)\right\},$ where $u_{n}$ is a solution to (I) of  $\eqref{eq1},$ is bounded in $W_{0}^{1, 2}(\Omega) .$
\end{lem}	\textbf{Proof.} It is sufficient to take $T_{k}\left(u_{n}\right)$ as a test function in (I) of problems \eqref{eq1}. $\square$
\begin{lem}\label{1}
	Let $0 \leq \theta <1,$ and let $f$ belong to $L^{m}(\Omega),$ $m \geq \max \left(\big(\frac{r+1}{1-\theta}\big)^{\prime}, \frac{2 N}{\theta(N-2)+N+2}\right)$ and $r>1,$ we have that:\\
	$\bullet$ the sequences $\left\{u_{n}\right\}$ and $\left\{\psi_{n}\right\}$ are bounded in $W_{0}^{1,2}(\Omega).$\\
	$\bullet$ the sequence $\left\{u_{n}\right\}$ is bounded in $L^{\sigma}(\Omega),$ with $\sigma=\max \left((1+\theta)m^{* *}, \frac{m(2 r+1+\theta)}{m+1} \right)$
	if $m<\frac{N}{2},$ and $\sigma=+\infty$ if $m>\frac{N}{2}.$
\end{lem}
\textbf{Proof.}
\textbf{$L^{\infty}(\Omega)$ estimate: }
Suppose that $m> \frac{N}{2}$, let $k>1$ and define $G_{k}(s)=(s-k)^{+}$. Choosing $G_{k}\left(u_{n}\right)$ as test function in \eqref{eq1}, we obtain, recalling \eqref{134},
$$
\begin{aligned}
\alpha \int_{\Omega}\left|\nabla G_{k}\left(u_{n}\right)\right|^{2} & \leq \int_{\Omega} M(x) \nabla G_{k}\left(u_{n}\right) \times \nabla G_{k}\left(u_{n}\right)
\end{aligned}
$$
\begin{equation}\label{333}
=\int_{\Omega} \frac{f_{n} G_{k}\left(u_{n}\right)}{(u_{n}+\frac{1}{n})^{\theta}} \leq \int_{\Omega} f G_{k}\left(u_{n}\right),
\end{equation}
where in the last passage we have used that $u_{n}+\frac{1}{n} \geq k \geq 1,$ on the set $\left\{u_{n} \geq k\right\}$ where $G_{k}\left(u_{n}\right) \neq 0 .$ Starting from inequality \eqref{333} and arguing as in \cite{4}, Théorème 4.2, we have that there exists a constant $C$ (independent on $n$ ), such that	
$$
\left\|u_{n}\right\|_{L^{\infty}(\Omega)} \leq C\|f\|_{L^{m}(\Omega)}.
$$
Since $u_{n}$ is bounded in $L^{\infty}(\Omega),$ as well.$\square$\\
\textbf{Estimates which use the lower order term In this step,} we will suppose that $m \geq \big(\frac{r+1}{1-\theta}\big)^{\prime} .$ Taking $u_{n}$ as test function in the first equation of \eqref{eq1}, using
\eqref{134} and dropping a positive term, we obtain
$$
\alpha \int_{\Omega}\left|\nabla u_{n}\right|^{2}+\int_{\Omega} \psi_{n}\left|u_{n}\right|^{r} \leq \int_{\Omega} f_{n} u_{n}^{1-\theta},
$$
while using $\psi_{n}$ as test function in (II) and \eqref{135}, we can see that
$$
\alpha \int_{\Omega}\left|\nabla \psi_{n}\right|^{2} \leq \int_{\Omega} \psi_{n}\left|u_{n}\right|^{r}.
$$
Thus we have, once again, that
\begin{equation}\label{233}
\alpha \int_{\Omega}\left|\nabla u_{n}\right|^{2}+\alpha \int_{\Omega}\left|\nabla \psi_{n}\right|^{2} \leq \int_{\Omega}\left|f_{n} \| u_{n}\right|^{1-\theta}.
\end{equation}
We now follow \cite{10}, let $\gamma \geq 1$ to be determined later, and choose $\left|u_{n}\right|^{2 \gamma-2} u_{n}$ as test function in the first equation of \eqref{eq1}$;$ using \eqref{134}, and dropping two
positive terms, we obtain, since $\left|f_{n}\right| \leq|f|$
\begin{equation}\label{244}
\alpha(2 \gamma-1) \int_{\Omega}\left|\nabla u_{n}\right|^{2}\left|u_{n}\right|^{2 \gamma-2} \leq \int_{\Omega}\left|f_{n}\right|\left|u_{n}\right|^{2 \gamma-1-\theta} \leq \int_{\Omega}\left|f \| u_{n}\right|^{2 \gamma-1-\theta}.
\end{equation}
On the other hand, taking $\left|u_{n}\right|^{\gamma}$ as a test function in (II), by estimate \eqref{135} and  using Young inequality, we obtain that
$$
\begin{aligned}
\int_{\Omega}\left|u_{n}\right|^{r+\gamma} &=\gamma \int_{\Omega} M(x) \nabla \psi_{n} \nabla u_{n}\left|u_{n}\right|^{\gamma-1}  \\
& \leq \beta \gamma \int_{\Omega}\left|\nabla \psi_{n}\left\|\nabla u_{n}\right\| u_{n}\right|^{\gamma-1}\\
&\leq C \int_{\Omega}\left|\nabla \psi_{n}\right|^{2}+C \int_{\Omega}\left|\nabla u_{n}\right|^{2}\left|u_{n}\right|^{2 \gamma-2}.
\end{aligned}
$$
Using \eqref{233} and \eqref{244} with this inequality, we deduce that
$$
\int_{\Omega}|u_{n}|^{r+\gamma} \leq C \int_{\Omega}|f||u_{n}|^{1-\theta}+C \int_{\Omega}| f|| u_{n}|^{2 \gamma-1-\theta},
$$
so that we have
\begin{equation}\label{255}
\int_{\Omega}\left|u_{n}\right|^{r+\gamma} \leq C \int_{\Omega}|f||u_{n}|^{1-\theta}+C \int_{\Omega}| f| |u_{n}|^{2 \gamma-1-\theta} ,
\end{equation}
where in the last passage, we have used that $2 \gamma-1-\theta \geq 1-\theta,$ since $\gamma \geq 1$. We now choose $\gamma=\frac{r(m-1)+m(\theta+1)}{m+1},$ so that $\gamma \geq 1$ since $m \geq \frac{r+1}{r+\theta}=(\frac{r+1}{1-\theta})^{\prime} .$ With this choice of $\gamma,$ we obtain $r+\gamma=\frac{m(2 r+1+\theta)}{m+1}=(2 \gamma-1-\theta) m^{\prime},$ so by Hölder inequality, we deduce from  \eqref{255} that
$$
\int_{\Omega}\left|u_{n}\right|^{\frac{m(2 r+1+\theta)}{m+1}} \leq C\|f\|_{L^{m}(\Omega)}\left[\int_{\Omega}\left|u_{n}\right|^{\frac{m(1-\theta)}{m-1}}\right]^{\frac{1}{m^{\prime}}}
+C\|f\|_{L^{m}(\Omega)}\left[\int_{\Omega}\left|u_{n}\right|^{\frac{m(2 r+1+\theta)}{m+1}}\right]^{\frac{1}{m^{\prime}}}.
$$
Thanks to the fact that $m>1,$ we, therefore, obtain (after simplifying equal terms) that:
$$
\left[\int_{\Omega}\left|u_{n}\right|^{\frac{m(2 r+1+\theta)}{m+1}}\right]^{\frac{1}{m}} \leq C\|f\|_{L^{m}(\Omega)},
$$
that is, the sequence $\left\{u_{n}\right\}$ is bounded in $L^{s}(\Omega),$ with $s=\frac{m(2 r+1+\theta)}{m+1} .$ As a consequence of this estimate, and of the fact that $s \geq m^{\prime},$ we have that
$$
\int_{\Omega}|f|\left|u_{n}\right|^{1-\theta} \leq C,
$$
so that from \eqref{233}, it follows that the sequences $\left\{u_{n}\right\}$ and $\left\{\psi_{n}\right\}$ are bounded in $W_{0}^{1,2}(\Omega).$ $\square$\\
\textbf{Estimates not using the lower order term in this step,} we will suppose that $m \geq \frac{2 N}{\theta(N-2)+N+2} .$ Let $u_{n}$ and $\psi_{n}$ be solutions of \eqref{eq1}, let $\gamma \geq 1,$ and take $\left|u_{n}\right|^{2 \gamma-2} u_{n}$ as test function in (I) of \eqref{eq1}, we have,  dropping two positive terms, and using \eqref{134},
$$
\alpha(2 \gamma-1) \int_{\Omega}\left|\nabla u_{n}\right|^{2}\left|u_{n}\right|^{2 \gamma-2} \leq \int_{\Omega} f_{n}\left|u_{n}\right|^{2 \gamma-2} u_{n}^{1-\theta}.
$$
By exploiting  Sobolev and Hölder inequalities, and since $\left|f_{n}\right| \leq|f|,$ we deduce
$$
\begin{aligned}
\frac{\alpha \mathcal{S}(2 \gamma-1)}{\gamma^{2}}\left[\int_{\Omega}\left|u_{n}\right|^{2^{*} \gamma} \right]^{\frac{2}{2^{*}}} & \leq \alpha(2 \gamma-1) \int_{\Omega}\left|\nabla u_{n}\right|^{2} u_{n}^{2 \gamma-2} \\
& \leq \int_{\Omega} f_{n}\left|u_{n}\right|^{2 \gamma-2} u_{n}^{1-\theta}\\ &\leq\|f\|_{L^{m}(\Omega)}\left[\int_{\Omega}\left|u_{n}\right|^{(2 \gamma-1-\theta) m^{\prime}}\right]^{\frac{1}{m^{\prime}}}.
\end{aligned}
$$
Imposing $ 2^{*} \gamma= (1+\theta)m^{* *},$ we have  $\gamma=\frac{(1+\theta)m^{* *}}{2^{*}},$ so that $\gamma \geq 1\left(\right.$ since $\left.(1+\theta)m^{* *} \geq 2^{*}\right)$ and $(2 \gamma-1-\theta) m^{\prime}=(1+\theta)m^{* *}=:s ,$ we have
$$
\left[\int_{\Omega}\left|u_{n}\right|^{s}\right]^{\frac{2}{2^{*}}} \leq C\|f\|_{L^{m}(\Omega)}\left[\int_{\Omega}\left|u_{n}\right|^{s}\right]^{\frac{1}{m^{\prime}}},
$$
so that
$$
\left[\int_{\Omega}\left|u_{n}\right|^{s}\right]^{\frac{1}{s}} \leq C\|f\|_{L^{m}(\Omega)}.
$$
Thus, the sequence $\left\{u_{n}\right\}$ is bounded in $L^{s}(\Omega),$
being $m \geq \frac{2 N}{N+2+\theta(N-2)},$ we have that the sequence $\left\{f_{n} u_{n}^{1-\theta}\right\}$ is bounded in $L^{1}(\Omega).$ Taking, $u_{n}$ as test function in the equation (I) of \eqref{eq1}, to obtain, after using  \eqref{134}
and dropping a positive term,
$$
\alpha \int_{\Omega}\left|\nabla u_{n}\right|^{2}+\int_{\Omega} \psi_{n}\left|u_{n}\right|^{r} \leq \int_{\Omega} f_{n} u_{n}^{1-\theta} \leq C,
$$
so that , the sequence $\left\{u_{n}\right\}$ is bounded in $W_{0}^{1,2}(\Omega),$ and the sequence $\left\{\psi_{n}\left|u_{n}\right|^{r}\right\}$ is bounded in $L^{1}(\Omega) .$ Choosing $\psi_{n}$ as test function in (II)  of \eqref{eq1}, and using \eqref{135}, we thus have:
$$
\alpha \int_{\Omega}\left|\nabla \psi_{n}\right|^{2} \leq \int_{\Omega} \psi_{n}\left|u_{n}\right|^{r} \leq C,
$$
so that also the sequence $\left\{\psi_{n}\right\}$ is bounded in $W_{0}^{1,2}(\Omega).$ $\square$\\
\subsection{ Proof of Theorem \ref{th21}}
In virtue of the Lemma \ref{1}, the sequence of approximated solutions $u_{n}$ is bounded in $ W_{0}^{1,2}(\Omega)\cap L^{\sigma}(\Omega)$. Therefore, there exists a function $u$ belongs to $ W_{0}^{1,2}(\Omega)\cap L^{\sigma}(\Omega)$ such that, up to subsequences, $u_{n}$ converges, weakly in $W_{0}^{1,2}(\Omega),$ weakly in $L^{\sigma}(\Omega),$ and almost everywhere in $\Omega,$ to some function $u,$ while $\psi_{n}$ converges, weakly in $W_{0}^{1,2}(\Omega)$ and almost everywhere in $\Omega,$ to some function $\psi$. Since the sequence $\left\{\left|u_{n}\right|^{r}\right\}$ is bounded in $L^{\rho}(\Omega),$ with $\rho=\frac{\sigma}{r}>1,$ it is weakly convergent in the same space to $|u|^{r}$. Therefore, one can pass to the limit in the identities
$$
\int_{\Omega} M(x) \nabla \psi_{n} \nabla w=\int_{\Omega}\left|u_{n}\right|^{r} w, \quad \forall w \in C_{c}^{1}(\Omega),
$$
to have that $\psi$ and $u$ are such that:
$$
\int_{\Omega} M(x) \nabla \psi \nabla w=\int_{\Omega}|u|^{r} w, \quad \forall w \in C_{c}^{1}(\Omega).
$$
Choosing $w=T_{k}(v),$ with $v \geq 0$ in $ C_{c}^{1}(\Omega),$ we arrive  at
$$
\int_{\Omega} M(x) \nabla \psi \nabla T_{k}(v)=\int_{\Omega}|u|^{r} T_{k}(v), \quad \forall k>0.
$$
Letting $ k$ tend to infinity, using Lebesgue theorem in the left-hand side (recall that $\psi$ belongs to $\left.W_{0}^{1,2}(\Omega)\right),$ and Beppo Levi theorem in the right-hand side, we deduce  that
$$
\int_{\Omega} M(x) \nabla \psi \nabla v=\int_{\Omega}|u|^{r} v, \quad \forall v \in C_{c}^{1}(\Omega), v \geq 0.
$$
If $v$ belongs to $C_{c}^{1}(\Omega),$ writing $v=v^{+}-v^{-},$ and subtracting the above identities written for $v^{+}$ and $v^{-}$ (not that both terms are finite, because the left-hand side is finite), we have that
$$
\int_{\Omega} M(x) \nabla \psi \nabla v=\int_{\Omega}|u|^{r} v, \quad \forall v \in C_{c}^{1}(\Omega),
$$
that is, $\psi$ is a weak solution of the second equation. We study now the first equation: We want to prove that
$\psi_{n}\left|u_{n}\right|^{r-1}$ strongly converges to $\psi|u|^{r-1}$ in $L^{1}(\Omega).$
First of all, let $\varepsilon>0, k>0,$ and choose $\frac{1}{\varepsilon} u_{n}^{+} T_{\varepsilon}\left(G_{k}\left(u_{n}\right)\right)$
as test function in the first equation of the system. Dropping two positive terms (those coming from the differential part of the equation), and using that $\left|f_{n}\right| \leq|f|,$ we obtain
$$
\begin{aligned}
\frac{1}{\varepsilon} \int_{\left\{u_{n} \geq k\right\}} \quad \psi_{n}\left[u_{n}^{+}\right]^{r}
T_{\varepsilon}\left(G_{k}\left(u_{n}\right)\right) \leq \frac{1}{\varepsilon} \int_{\left\{u_{n} \geq k\right\}}|f_{n}||u_{n}|^{1-\theta}
T_{\varepsilon}\left(G_{k}\left(u_{n}\right)\right) \\
\quad \leq \frac{1}{\varepsilon} \int_{\left\{u_{n} \geq k\right\}}|f||u_{n}|^{1-\theta} T_{\varepsilon}\left(G_{k}\left(u_{n}\right)\right).
\end{aligned}
$$
Letting $\varepsilon$ tend to zero, using Fatou lemma on the left-hand side, and Lebesgue theorem on the right-hand one (recall that every $u_{n}$ is a function in $\left.L^{\infty}(\Omega)\right)$, we have that
$$
\begin{aligned}
\int_{\left\{u_{n} \geq k\right\}} \psi_{n}\left[u_{n}^{+}\right]^{r} & \leq \int_{\left\{u_{n} \geq k\right\}}
|f|u_{n} \leq\left[\int_{\left\{u_{n} \geq k\right\}}|f|^{m}\right]^{\frac{1}{m}}
||u_{n}^{1-\theta}||_{L^{m^{\prime}}(\Omega)}.\\
& \leq C\left[\int_{\left\{u_{n} \geq k\right\}}|f|^{m}\right]^{\frac{1}{m}},
\end{aligned}
$$
since the sequence $\left\{u_{n}\right\}$ is bounded in $L^{m^{\prime}}(\Omega)$ being $\sigma \geq (1-\theta) m^{\prime}.$ Then
$$
\int_{\left\{u_{n} \geq k\right\}} \psi_{n}\left|u_{n}\right|^{r} \leq C\left[\int_{\left\{u_{n} \geq k\right\}}|f|^{m}\right]^{\frac{1}{m}}.
$$
Let now $E$ be a measurable subset of $\Omega$. So that
$$
\begin{aligned}
\int_{E} \psi_{n}u_{n}^{r} &=\int_{E \cap\left\{u_{n} \leq k\right\}} \psi_{n}\left|u_{n}\right|^{r}+\int_{E \cap\left\{u_{n} \geq k\right\}} \psi_{n}u_{n}^{r} \\
& \leq k^{r} \int_{E} \psi_{n}+C\left[\int_{\left\{u_{n} \geq k\right\}}|f|^{m}\right]^{\frac{1}{m}}.
\end{aligned}
$$
Now we choose $\varepsilon>0$, and let $k$ large enough,we obtain
$$
C\left[\int_{\left\{u_{n} \geq k\right\}}|f|^{m}\right]^{\frac{1}{m}} \leq \varepsilon, \quad \forall n \in \mathbb{N}.
$$
Such a choice of $k$ is possible, since the measure of $\left\{u_{n} \geq k\right\}$ tends to zero as $k$ tends to infinity, uniformly in $n,$ as a consequence of the boundedness of $\left\{u_{n}\right\}$ in (for example) $L^{1}(\Omega),$ and since $|f|^{m}$ belongs to $L^{1}(\Omega) .$ Once $k$ has been chosen, let $\delta>0$ be such that meas $(E) \leq \delta$ implies that:
$$
k^{r} \int_{E} \psi_{n} \leq \varepsilon, \quad \forall n \in \mathbb{N}.
$$
Such a choice of $\delta$ is possible thanks to Vitali theorem, since the sequence $\left\{\psi_{n}\right\}$ is strongly convergent in (at least) $L^{1}(\Omega)$ being bounded in $W_{0}^{1,2}(\Omega) .$ Thus, the sequence $\left\{\psi_{n}\left|u_{n}\right|^{r}\right\}$ is uniformly equi-integrable. Since it is almost everywhere convergent, Vitali theorem implies that:
$$\psi_{n}\left|u_{n}\right|^{r} \mbox{ strongly converges to} \,\,\psi|u|^{r}\,\, \mbox{in}\,\, L^{1}(\Omega).$$
With the same technique, one can prove that the sequence $\left\{\psi_{n}\left|u_{n}\right|^{r-1}\right\}$ is uniformly equi-integrable, so that $\psi_{n}\left|u_{n}\right|^{r-1}$ strongly converges to $\psi|u|^{r-1}$ in $L^{1}(\Omega).$
We want  to pass to the limit in (I) of \eqref{eq1}.
For the limit of the right hand of (I) in \eqref{eq1}. Let $w= \{\varphi \neq 0\}$ then
by Lemma \ref{lem22} (see the appendix ), one has, for $\varphi$ in $C_{c}^{1}(\Omega),$ we have that
$$
0 \leq\left|\frac{f_{n} \varphi}{(u_{n}+\frac{1}{n})^{\theta}}\right| \leq \frac{\|\varphi\|_{L^{\infty}(\Omega)}}{c_{\omega}^{\theta}} f.
$$
Therefore, by Lebesgue convergence  Theorem, we obtain
$$
\lim _{n \rightarrow+\infty} \int_{\Omega} \frac{f_{n} \varphi}{(u_{n}+\frac{1}{n})^{\theta}}=\int_{\Omega} \frac{f \varphi}{u^{\theta}}.
$$
On other hand, by Lemma \ref{lc1}, we deduce $ T_{k}(u_{n}) \rightharpoonup T_{k}(u)$ weakly in
$W_{0}^{1,2}(\Omega).$ Then by Proposition 4.1 in \cite{c1} and Theorem 2.3 in \cite{c11}, we obtain $ \nabla u_{n}$ converges to $\nabla u$
almost everywhere in $\Omega.$
Now, we can pass to the limit in the identities:
$$
\int_{\Omega}a(x) \nabla u_{n} \nabla \eta+\int_{\Omega} \psi_{n}\left|u_{n}\right|^{r-2} u_{n} \eta=\int_{\Omega} \frac{f_{n}}{(\frac{1}{n}+ u_{n})^{\theta}} \eta, \quad \forall \eta \in  C_{c}^{1}(\Omega),
$$
to have that
$$
\int_{\Omega}a(x) \nabla u \nabla \eta+\int_{\Omega} \psi|u|^{r-2} u \eta=\int_{\Omega} \frac{f }{u^{\theta}}\eta, \quad \forall \eta \in C_{c}^{1}(\Omega),
$$
as desired.
\section{Saddle points}
In this section, we can prove that the solution $(u, \psi)$ of system
\eqref{133} given by Theorem \ref{th21} can be seen (under some assumptions on $r$ and $f)$ as a saddle point of a suitable functional.
\begin{rem}\label{22} If $1<r \leq \frac{N+2+(N-2)\theta}{N-2},$ and $f$ belongs to $L^{m}(\Omega),$ with $m \geq (\frac{2^{*}}{1-\theta})^{'},$ then not only $\psi$ but also $u$ is a weak solution of the first equation of $\eqref{eq1} .$ Indeed, since both $u$ and $\psi$ belong to $L^{2^{*}}(\Omega)$ (being $W_{0}^{1,2}(\Omega)$ functions), we have that:
	$$
	\psi|u|^{r-2} u \in L^{\rho}(\Omega), \quad \rho=\frac{2^{*}}{r}
	$$
	since, by the assumptions on $r$
	$$
	\frac{2^{*}}{r} \geq \frac{2 N}{N-2} \frac{N-2}{N+2+(N-2)\theta}=(\frac{2^{*}}{1-\theta})^{'},
	$$
	the function $\psi|u|^{r-2} u$ belongs to the dual of $W_{0}^{1,2}(\Omega) ;$ therefore, one has (by density of $W_{0}^{1,2}(\Omega) \cap L^{\infty}(\Omega)$ in $\left.W_{0}^{1,2}(\Omega)\right)$
	$$
	\int_{\Omega}a(x) \nabla u \nabla \varphi+\int_{\Omega} \psi|u|^{r-2} u \varphi=\int_{\Omega} \frac{f }{u^{\theta}}\varphi, \quad \forall \varphi \in C_{c}^{1}(\Omega),
	$$
	as desired.
\end{rem}
Thanks to this remark, we have the following  theorem:
\begin{thm}
	Suppose that a and $M$ satisfy \eqref{134} and $\eqref{135},$ and that $M$ is symmetric. Let $1<r \leq \frac{2 N}{N+2+(N-2)\theta}$ and let $f$ in $L^{m}(\Omega),$ with $m \geq (\frac{2^{*}}{1-\theta})^{'}.$ Then, the solution $(u,\psi)$ of system \eqref{133} given by Theorem \ref{th21} is a saddle point of the functional $J$ defined in \eqref{1.6};  that is
	\begin{equation}\label{2100}
	J(u, \varphi) \leq J(u, \psi) \leq J(v, \psi), \quad \forall v, \varphi \in W_{0}^{1,2}(\Omega) \text { such that } \psi|v|^{r} \in L^{1}(\Omega).
	\end{equation}
\end{thm}
\textbf{Proof.} We begin with the second equation of \eqref{133}$;$ by Theorem $\ref{th21}, \psi$ is a weak solution of the second equation of $\eqref{133} .$ Choosing $\frac{\psi-\varphi^{+}}{r},$ with $\varphi$ in $W_{0}^{1,2}(\Omega),$ as test function, we get
$$
\frac{1}{r} \int_{\Omega} M(x) \nabla \psi \nabla\left(\psi-\varphi^{+}\right)=\frac{1}{r} \int_{\Omega}|u|^{r}\left(\psi-\varphi^{+}\right).
$$
Adding and subtracting the term
$$
\frac{1}{2 r} \int_{\Omega} M(x) \nabla \varphi^{+} \nabla \varphi^{+},
$$
we have, after straightforward passages
$$
\begin{aligned}
\frac{1}{2 r} \int_{\Omega} M(x) \nabla\left(\psi-\varphi^{+}\right) \nabla\left(\psi-\varphi^{+}\right)+\frac{1}{2 r} \int_{\Omega} M(x) \nabla \psi \nabla \psi-\frac{1}{r} \int_{\Omega} \psi|u|^{r} \\
\quad=\frac{1}{2 r} \int_{\Omega} M(x) \nabla \varphi^{+} \nabla \varphi^{+}-\frac{1}{r} \int_{\Omega} \varphi^{+}|u|^{r}
\end{aligned}
$$
since the first term is positive, we, therefore, have that (recall that $\psi \geq 0,$ so that $\left.\psi=\psi^{+}\right)$
$$
\frac{1}{2 r} \int_{\Omega} M(x) \nabla \psi \nabla \psi-\frac{1}{r} \int_{\Omega} \psi^{+}|u|^{r} \leq \frac{1}{2 r} \int_{\Omega} M(x) \nabla \varphi^{+} \nabla \varphi^{+}-\frac{1}{r} \int_{\Omega} \varphi^{+}|u|^{r},
$$
for every $\varphi$ in $W_{0}^{1,2}(\Omega) .$ Changing sign to this identity, and adding to both sides the (finite, thanks to the assumptions on $f$ and to the fact that $u$ belongs to $\left.W_{0}^{1,2}(\Omega)\right)$ term
$$
\frac{1}{2} \int_{\Omega} a(x)|\nabla u|^{2}-\frac{1}{1-\theta}\int_{\Omega} f u^{1-\theta},
$$
we arrive
$$
J(u, \varphi) \leq J(u, \psi), \quad \forall \varphi \in W_{0}^{1,2}(\Omega),
$$
which is the first half of $\eqref{2100} .$ As for the second, by Remark \ref{22}, we obtain that $u$ is a weak solution of the first equation of $\eqref{133}.$ Fix $\psi \in W_{0}^{1, 2}(\Omega)$ and let $I$ be the functional defined on $W_{0}^{1, 2}(\Omega)$ as $I(v):=J(v, \psi)$.
If the matrix $M(x)$ and a(x) is symmetric, and if $f$ belongs to $L^{m}(\Omega),$ with $m>\left(\frac{2^{*}}{1-\theta}\right)^{\prime}$ the solution of \eqref{133} given by Theorem \ref{th21} is the minimum of the functional
$$
I(v)=\frac{1}{2} \int_{\Omega} a(x) \nabla v \times \nabla v-\frac{1}{2 r}\int_{\Omega} M(x) \nabla \psi \nabla \psi
$$
$$
+\frac{1}{r} \int_{\Omega} \psi^{+}|v|^{r} -\frac{1}{1-\theta} \int_{\Omega} f v^{1-\theta}, \quad v \in W_{0}^{1,2}(\Omega)
$$
which is well defined since $\theta<1$. Indeed, if we consider the functional
$$
I_{n}(v)=\frac{1}{2} \int_{\Omega} a(x) \nabla v \times \nabla v - \int_{\Omega} M(x) \nabla \psi \nabla \psi$$
$$+\frac{1}{r} \int_{\Omega} \psi^{+}|v|^{r} -\frac{1}{1-\theta} \int_{\Omega} f_{n}\left(v^{+}+\frac{1}{n}\right)^{1-\theta}, \quad v \in W_{0}^{1,2}(\Omega)
$$
with $f_{n}=\min (f(x), n),$ then there exists a minimum $u_{n}$ of $I_{n} .$ From the inequality $I_{n}\left(u_{n}\right) \leq$ $I_{n}\left(u_{n}^{+}\right)$ one can prove that $u_{n} \geq 0,$ so that $u_{n}$ is a solution of the Euler equation for $I_{n},$ i.e., of \eqref{133}. Therefore, by Lemma \ref{lem22} and Remark \ref{r23}$, u_{n}$ is unique and increasing in $n,$ satisfies \eqref{223} and, from the inequality $I\left(u_{n}\right) \leq I_{n}(0) \leq C,$ it is bounded in $W_{0}^{1,2}(\Omega)$ (with the same proof of  Lemma \ref{1} ). If $u$ is the limit of $u_{n}$, letting $n$ tend to infinity in the inequalities $I_{n}\left(u_{n}\right) \leq I_{n}(v),$ one finds that $I(u) \leq I(v),$ so that $u$ is a minimum of $I,$ and $u$ is a solution of \eqref{133} (by Theorem \ref{th21} ). Since $u$ satisfies $\eqref{223} ,$ Eq. \eqref{133} can be seen as the Euler equation for $ I ;$ note that $I$ is not differentiable on $W_{0}^{1,2}(\Omega) .$ We obtain that:
$$
J(u, \psi) \leq J(v, \psi), \quad \forall v \in W_{0}^{1,2}(\Omega) \text { such that } \psi|v|^{r} \in L^{1}(\Omega),
$$
which is the second part of \eqref{2100}.
\section{Appendix: Basic Results and Existence for Bounded Data }
In this Appendix, we will prove some results concerning the first equation of
system \eqref{133}, and the whole system in the case of bounded data.\\
Now we prove the existence of a solution to the following approximating
problems:
\begin{equation}\label{P3}
\left\{\begin{array}{ll}
-\operatorname{div}(a(x) u_{n})+g(x)\left|u_{n}\right|^{r-2} u_{n}=\frac{f_{n}(x)}{\left(\left|u_{n}\right|+\frac{1}{n}\right)^{\theta}} & \text { in } \Omega \\
u_{n}>0 & \text { in }  \Omega\\
u_{n}=0 & \text { on } \partial \Omega
\end{array}\right.
\end{equation}

where $\Omega$ is a bounded open subset of $\mathbb{R}^{N}, N \geq 2, f$ is a positive (that is $f(x) \geq 0$ and not zero a.e.) function in $L^{m}(\Omega),$ with $m \geq 1,$  $0< \theta<1$ and
$g(x)\in L^{1}(\Omega),$ with
\begin{equation}
0 < \lambda \leq g(x).
\end{equation}
Due to the nature of the approximation, the sequence un will be increasing with n,
so that the (strict) positivity of the limit will be derived from the (strict) positivity
of any of the un (which in turn will follow by the standard maximum principle for
elliptic equations).
\begin{lem}\label{lem1}
	Problem \eqref{P3} has a nonnegative solution $u_{n}$ in $W_{0}^{1,2}(\Omega)\cap L^{\infty}(\Omega)$.
\end{lem}
In order to prove Lemma \ref{lem1} , we will work by approximation, namely by introducing the following
\begin{equation}\label{p22}
\left\{\begin{array}{ll}
-\operatorname{div} (a(x) u_{n, k})+g(x)T_{k}\left(\left|u_{n, k}\right|^{r-2} u_{n, k}\right)=\frac{f_{n}(x)}{\left(\left|u_{n, k}\right|+\frac{1}{n}\right)^{\theta}} & \text { in } \Omega \\
u_{n, k}=0 & \text { on } \Omega
\end{array}\right.
\end{equation}
where $n, k \in \mathbb{N}, 0 \leq f_{n}(x):=T_{n}(f(x)) \in L^{\infty}(\Omega), 0< \theta < 1$ and $r \geq 1$.
Thanks to [\cite{le} $, Thoerem \,2]$, we know that there exists $u_{n, k} \in W_{0}^{1,2}(\Omega)$ weak solution to \eqref{p22} for each $n, k \in \mathbb{N}$ fixed. Moreover $u_{n, k} \in L^{\infty}(\Omega)$ for all $n, k \in \mathbb{N}$ since, if $m \geq 1$ is fixed, taking
$G_{m}\left(u_{n, k}\right) \in W_{0}^{1,2}(\Omega)$ as test function in \eqref{p22} and using that $G_{m}\left(u_{n, k}\right)$ and $T_{k}\left(\left|u_{n, k}\right|^{r-2} u_{n, k}\right)$
have the same sign of $u_{n, k}$, we immediately find that
$$
\alpha\int_{\Omega}\left|\nabla G_{m}\left(u_{n, k}\right)\right|^{2} \leq \int_{\Omega} f_{n} G_{m}\left(u_{n, k}\right)
$$
and so we can proceed as in \cite{p7} to end up with $u_{n, k} \in L^{\infty}(\Omega).$ Moreover the previous $L^{\infty}$ estimate is independent from $k \in \mathbb{N}$. Now taking $u_{n, k}$ as a test function in the weak formulation of \eqref{p22}, we find that $u_{n, k}$ is bounded in $W_{0}^{1,2}(\Omega)$ with respect to $k$ for $n \in \mathbb{N}$ fixed. Since $u_{n, k}$ is bounded in $L^{\infty}(\Omega)$ independently on $k$, for each $n \in \mathbb{N}$ fixed we choose $k_{n}$ large enough to obtain the following scheme of approximation
\begin{equation}\label{p23}
\left\{\begin{array}{ll}
-\operatorname{div}(a(x) u_{n})+g(x)\left|u_{n}\right|^{r-2} u_{n}=\frac{f_{n}(x)}{\left(\left|u_{n}\right|+\frac{1}{n}\right)^{\theta}} & \text { in } \Omega \\
u_{n}=0 & \text { on } \partial \Omega
\end{array}\right.
\end{equation}
where $u_{n} \in W_{0}^{1,2}(\Omega) \cap L^{\infty}(\Omega)$ is given by $u_{n, k_{n}}$.
As concerns the sign of $u_{n}$, taking $u_{n}^{-}:=\min \left(u_{n}, 0\right) \in W_{0}^{1,2}(\Omega) \cap L^{\infty}(\Omega)$ as test function in \eqref{p23}, we find
$$
\int_{\Omega}a(x)\left|\nabla u_{n}^{-}\right|^{2}+\int_{\Omega}g(x)\left|u_{n}\right|^{r-2}\left(u_{n}^{-}\right)^{2}=\int_{\Omega} \frac{f_{n}}{\left(\left|u_{n}\right|+\frac{1}{n}\right)^{\theta}} u_{n}^{-} \leq 0
$$
and so that $u_{n} \geq 0$ almost everywhere in $\Omega$.
\begin{lem}\label{lem22}
	The sequence $u_{n}$ is increasing with respect to $n, u_{n}>0$ in $\Omega$, and for every $\omega \subset \subset \Omega$ there exists $c_{\omega}>0$ (independent on $n$ ) such that
	\begin{equation}\label{223}
	u_{n}(x) \geq c_{\omega}>0 \quad \mbox{ for every }\,\, x\,\, \mbox{ in} \,\, \omega,\mbox{ for every}\,\, n\,\,  \mbox{in}\,\,  \mathbb{N}.
	\end{equation}
	Moreover there exists the pointwise limit $u \geq c_{\omega}$ of the sequence $u_{n} .$\\
\end{lem}
\textbf{Proof.} Since $0 \leq f_{n} \leq f_{n+1}$ and $\theta>0$, one has (distributionally)
$$
-\operatorname{div}\left(a(x) \nabla u_{n}\right)+g(x)\left|u_{n}\right|^{r-2} u_{n}=\frac{f_{n}}{\left(u_{n}+\frac{1}{n}\right)^{\theta}} \leq \frac{f_{n+1}}{\left(u_{n}+\frac{1}{n+1}\right)^{\theta}},
$$
so that
$$
\begin{array}{l}
-\operatorname{div}\left(a(x)\left( \nabla u_{n} - \nabla u_{n+1}\right)\right)  +g(x)\left( \left|u_{n}\right|^{r-2} u_{n}-\left|u_{n+1}\right|^{r-2} u_{n+1} \right) \\
\quad \leq f_{n+1} \frac{\left(u_{n+1}+\frac{1}{n+1}\right)^{\theta}-\left(u_{n}+\frac{1}{n+1}\right)^{\theta}}{\left(u_{n}+\frac{1}{n+1}\right)^{\theta}\left(u_{n+1}+\frac{1}{n+1}\right)^{\theta}}.
\end{array}
$$
We now choose $\left(u_{n}-u_{n+1}\right)^{+}$ as test function and taking into  account  the monotonicity of the function $t \rightarrow|t|^{r-2} t .$ For the right hand side we observe that
$$
\left[\left(u_{n+1}+\frac{1}{n+1}\right)^{\theta}-\left(u_{n}+\frac{1}{n+1}\right)^{\theta}\right]\left(u_{n}-u_{n+1}\right)^{+} \leq 0,
$$
recalling that $f_{n+1} \geq 0$, we obtain
$$
0 \leq \alpha \int_{\Omega}\left|\nabla\left(u_{n}-u_{n+1}\right)^{+}\right|^{2} \leq 0.
$$
Therefore $\left(u_{n}-u_{n+1}\right)^{+}=0$ almost everywhere in $\Omega$, which implies $u_{n} \leq u_{n+1}$. Since $u_{1}$ belongs to $L^{\infty}(\Omega)$, and there exists a constant (only depending on $\Omega$ and $N$ ) such that
$$
\left\|u_{1}\right\|_{L^{\infty}(\Omega)} \leq C\left\|f_{1}\right\|_{L^{\infty}(\Omega)} \leq C
$$
one has
$$
\begin{array}{l}
-\operatorname{div}\left(a(x)\nabla u_{1}\right)+g(x)\left|u_{1}\right|^{r-2} u_{1} \\
\quad=\frac{f_{1}}{\left(u_{1}+1\right)^{\theta}} \geq \frac{f_{1}}{\left(\left\|u_{1}\right\|_{L^{\infty}(\Omega)}+1\right)^{\theta}} \geq \frac{f_{1}}{(C+1)^{\theta}}.
\end{array}
$$
Since $\frac{f_{1}}{(C+1)^{\theta}}$ is not identically zero, the strong maximum principle implies that $u_{1}>0$ in $\Omega$ (see \cite{p17}; observe that $u_{1}$ is differentiable by Chapter 4 of  \cite{p12} , and that \eqref{223} holds for $u_{1}$ (with $c_{\omega}$ only depending on $\omega, N, f_{1}$ and $\theta$ ). Since $u_{n} \geq u_{1}$ for every $n$ in $\mathbb{N},$\eqref{223} holds for $u_{n}$ (with the same constant $c_{\omega}$ which is then independent on $n$ ). 	
\begin{rem}\label{r23}
	If $u_{n}$ and $v_{n}$ are two solutions of \eqref{p23}, repeating the argument of the first part of the proof of Lemma \ref{lem22} shows that $u_{n} \leq v_{n} .$ By symmetry, this implies that the solution of \eqref{p23} is unique.
\end{rem}

\begin{thm}\label{3.4}
	Let $n\in \mathbb{N},$ $f$ be a function in $L^{\infty}(\Omega),$ and let $r>1 .$ Then, there exist $u$ and $\varphi,$ weak solutions of the system
	\begin{equation}\label{ii1}
	\left\{\begin{array}{l}
	u_{n} \in W_{0}^{1,2}(\Omega):-\operatorname{div}\left(a(x) \nabla u_{n}\right)+\varphi_{n}|u_{n}|^{r-2} u_{n}=\frac{f}{(\frac{1}{n}+u_{n})^{\theta}}, \\
	\varphi_{n} \in W_{0}^{1,2}(\Omega):-\operatorname{div}(M(x) \nabla \varphi_{n})=|u_{n}|^{r}
	\end{array}\right.
	\end{equation}
	Furthermore, $u_{n}$ and $\varphi_{n}$ belong to $L^{\infty}(\Omega),$ $\varphi_{n} > 0,$ 	$u_{n}>0 $ and $ 0< \theta <1.$
\end{thm}
\textbf{Proof.}
Fix $ \psi_{n} \in  W_{0}^{1, 2}(\Omega),$ let $n\in \mathbf{N}$  and we define $S: W_{0}^{1, 2}(\Omega) \rightarrow W_{0}^{1, 2}(\Omega)$ as the operator such that $v_{n}=S(\psi_{n})$.  By the maximum principle, $  \psi_{n} > 0 $, taking account Lemma \ref{lem1}  and Remark \ref{r23}, there, exists a unique solution $v_{n}$
of:
\begin{equation}\label{23}
-\operatorname{div}\left(a(x) \nabla v_{n}\right)+\psi_{n}|v_{n}|^{r-2} v=\frac{f}{(\frac{1}{n}+v_{n})^{\theta}}.
\end{equation}
Since, by Lemma \ref{lem1}, one has
\begin{equation}\label{24}
\|v_{n}\|_{W_{0}^{1,2}(\Omega)} \leq C_{1}\|f\|_{L^{\infty}(\Omega)}, \quad\|v_{n}\|_{L^{\infty}(\Omega)} \leq C_{1}\|f\|_{L^{\infty}(\Omega)} .
\end{equation}
Now we define $T: W_{0}^{1, 2}(\Omega) \rightarrow W_{0}^{1, 2}(\Omega)$ as the operator such that $\zeta_{n}=T(v_{n})=T(S(\psi_{n}))$.
Thanks to the results in \cite{03}, $\zeta_{n}$ is the unique weak solution of the Euler-Lagrange equation
\begin{equation}\label{25}
-\operatorname{div}\left(M(x)|\nabla \zeta_{n}|\right)=|v_{n}|^{r} , \quad \zeta_{n} \in W_{0}^{1, 2}(\Omega)
\end{equation}
Following \cite{02}, we thus have
$$
\|\zeta_{n}\|_{W_{0}^{1, 2}(\Omega)}+\|\zeta_{n}\|_{L^{\infty}(\Omega)} \leq C_{2}\|v_{n}\|_{L^{\infty}(\Omega)}^{r},
$$
using \eqref{24}, we deduce that,
\begin{equation}\label{27}
\|\zeta_{n}\|_{W_{0}^{1, 2}(\Omega)}+\|\zeta_{n}\|_{L^{\infty}(\Omega)} \leq C\|f\|_{L^{\infty}(\Omega)}=: R
\end{equation}
where $C_{1}$ and $C_{2}$ are positive constants not depending on $v_{n} .$

We want to prove that $T \circ S$ has a fixed point by Schauder's fixed point theorem. By \eqref{27} we have that $\overline{B_{R}(0)} \subset W_{0}^{1, 2}(\Omega)$ is invariant for $T \circ S .$ Let $ \psi_{k}=:(\psi_{n,k})_{k} \subset W_{0}^{1, 2}(\Omega)$ be a sequence weakly convergent to some $\psi$ and let $v_{k}=:(v_{n,k})_{k}=S\left(\psi_{k}\right) .$ As a consequence of \eqref{24}, there exists a subsequence indexed by $v_{k}$ \,\,\,\mbox{such that}\,\,
\begin{equation}\label{28}
v_{k} \rightarrow v\,\,\,\mbox{weakly in}\,\, W_{0}^{1, 2}(\Omega),  \,\,\,\mbox{and a.e. in}\,\, \Omega
\end{equation}
$$
v_{k} \rightarrow v \text { weakly-* in } L^{\infty}(\Omega).
$$
Moreover, we have
$$
-\operatorname{div}\left( a(x)\nabla v_{k}\right)=\frac{f}{(\frac{1}{n}+v_{k})^{\theta} }-\left(\psi_{k}^{+}\right)\left|v_{k}\right|^{r-2} v_{k}=: g_{k}
$$
and, using Hölder's inequality, the Poincaré inequality and \eqref{24}, we obtain

$$
\left\|g_{k}\right\|_{L^{1}(\Omega)} \leq C \|f\|_{L^{\infty}(\Omega)}+\left\|v_{k}\right\|_{L^{\infty}(\Omega)}^{r-1}\left\|\psi_{k}\right\|_{L^{1}(\Omega)}$$
$$ \leq C\|f\|_{L^{\infty}(\Omega)}+ C_{1}\|f\|_{L^{\infty}(\Omega)}^{r-1}\left\|\psi_{k}\right\|_{W_{0}^{1, 2}(\Omega)} \leq C .
$$
Then, by Theorem 2.1 in \cite{01}, we obtain that $\nabla v_{k}$ converges to $\nabla v_{n}$ almost everywhere in $\Omega$. Since
$$
\left\|\nabla v_{k}\right\|_{\left(L^{2}(\Omega)\right)^{N}}=\left\|v_{k}\right\|_{W_{0}^{1, 2}(\Omega)} \leq C_{1}\|f\|_{L^{m}(\Omega)},
$$
thus, we conclude that
\begin{equation}\label{29}
\nabla v_{k} \rightarrow \nabla v_{n} \text { weakly in }\left(L^{2}(\Omega)\right)^{N} \text { . }
\end{equation}
We recall that $v_{k}$ satisfies
$$
\int_{\Omega}a(x) \nabla v_{k} \cdot \nabla w+ \int_{\Omega}\psi_{k} \left|v_{k}\right|^{r-2} v_{k} w=\int_{\Omega} \frac{f}{(\frac{1}{n}+v_{k})^{\theta}} w, \quad \forall w\in C_{c}^{1}(\Omega).
$$
Letting $k$ tend to infinity, by \eqref{28},\eqref{29} and Vitali's theorem, we have that
$$
\int_{\Omega}| \nabla v_{n} \cdot \nabla w+\int_{\Omega}\psi_{n}|v_{n}|^{r-2} v_{n} w=\int_{\Omega}\frac{ f}{(\frac{1}{n}+v_{n})^{\theta}} w, \quad \forall w \in C_{c}^{1}(\Omega),
$$
so that $v$ is the unique weak solution of \eqref{23} and it does not depend on the subsequence. Hence $v_{k}=S\left(\psi_{k}\right)$ converges to $v_{n}=S(\psi_{n})$ weakly in $W_{0}^{1, 2}(\Omega)$ and weakly-* in $L^{\infty}(\Omega) .$ Then
\begin{equation}\label{210}
\left|v_{k}\right|^{r} \rightarrow|v_{n}|^{r} \text { strongly in } L^{q}(\Omega)\,\,\, \forall q<+\infty \text { and }\left\|\left|v_{k}\right|^{r} \right\|_{L^{1}(\Omega)} \leq C
\end{equation}
Thanks to \eqref{27},\eqref{210} and proceeding in the same way, we get
\begin{equation}\label{211}
\zeta_{k}:=\zeta_{n,k} =T\left(v_{k}\right) \rightarrow \zeta_{n}=T(v_{n})\,\,\mbox{ weakly in }\,W_{0}^{1, 2}(\Omega),\,\, \mbox{and weakly-* in }\,L^{\infty}(\Omega)
\end{equation}
$$
\left|\nabla \zeta_{k}\right| \nabla \zeta_{k} \rightarrow|\nabla \zeta_{n}| \nabla \zeta_{n} \text { weakly in }\left(L^{2}(\Omega)\right)^{N}.
$$
and $\zeta$ is the unique weak solution of \eqref{25} . Now we want to prove that $\zeta_{k}$ converges to $\zeta$ strongly in $W_{0}^{1, 2}(\Omega) .$ In order to obtain this, by Lemma 5 in \cite{05} , it is sufficient to prove the following
\begin{equation}\label{212}
\lim _{k \rightarrow \infty} \int_{\Omega }| \nabla\left(\zeta_{k}-\zeta_{n}\right)|^{2}=0.
\end{equation}
We have that
\begin{equation}\label{213}
\begin{aligned}
\int_{\Omega}\left(\left|\nabla \zeta_{k}\right| -|\nabla \zeta_{n}| \right) \cdot \nabla\left(\zeta_{k}-\zeta_{n}\right) &=\int_{\Omega}\left|\nabla \zeta_{k}\right|^{2}-\int_{\Omega}|\nabla \zeta_{n}|  \cdot \nabla \zeta_{k} \\
&-\int_{\Omega}\left|\nabla \zeta_{k}\right|  \cdot \nabla \zeta_{n}+\|\zeta_{n}\|_{W_{0}^{1, 2}(\Omega)}^{2}
\end{aligned}
\end{equation}
The second and the third term on the right hand side of \eqref{213} converge, by \eqref{211}, to $\|\zeta_{n}\|_{W_{0}^{1, 2}(\Omega)}^{2} .$ Then it is sufficient to prove that
\begin{equation}\label{214}
\lim _{k \rightarrow \infty}\left\|\zeta_{k}\right\|_{W_{0}^{1, 2}(\Omega)}^{2}=\|\zeta_{n}\|_{W_{0}^{1, 2}(\Omega)}^{2} .
\end{equation}
Since $\zeta_{k}$ is equal to $T\left(v_{k}\right) \geq 0$, we deduce that
$$
\int_{\Omega}\left|\nabla \zeta_{k}\right|^{2}=\int_{\Omega}\left|v_{k}\right|^{r} \zeta_{k}.
$$
Using  Vitali's Theorem and \eqref{210}, we have that
$$
\lim _{k \rightarrow \infty} \int_{\Omega}\left|v_{k}\right|^{r} \zeta_{k}=\int_{\Omega}|v_{n}|^{r} \zeta=\|\zeta_{n}\|_{W_{0}^{1, 2}(\Omega)}^{2},$$
so that \eqref{214} is true and \eqref{212} is proved. Hence we have proved that if $\psi_{k}$ converges to $\psi_{n}$ weakly in $W_{0}^{1, 2}(\Omega)$ then $\zeta_{k}=T\left(S\left(\psi_{k}\right)\right)$ converges to $\zeta_{n}=T(S(\psi_{n}))$ strongly in $W_{0}^{1, 2}(\Omega) .$ As a consequence we have that $T \circ S$ is a continuous operator and that $T\left(S\left(\overline{B_{R}(0)}\right)\right) \subset W_{0}^{1, 2}(\Omega)$ is a compact subset. Then there exists, by Schauder's fixed point Theorem, a function $\varphi_{n}$ in $W_{0}^{1, 2}(\Omega)$ such that $\varphi_{n}=T(S(\varphi_{n}))$ and, since $T(v_{n}) \geq 0$ for every $v_{n}$ in $W_{0}^{1, 2}(\Omega), \varphi_{n}$ is nonnegative. Moreover let $u_{n}=S(\varphi_{n})$, we have that $u_{n}$ is a weak solution of \eqref{ii1} .

\end{document}